\UseRawInputEncoding
\documentclass[12pt]{article}
\usepackage{}
\usepackage{mathrsfs}
\usepackage{amssymb}
\usepackage{amsmath}
\usepackage{amsfonts}
\usepackage{amstext}
\usepackage{amsthm}
\usepackage{graphicx}
\usepackage{subfig}
\usepackage{color}
\usepackage{indentfirst,latexsym,bm}
\usepackage{mathrsfs}
\usepackage{cite}
\usepackage[all]{xy}
\usepackage{newlfont}

\usepackage{enumerate}
\numberwithin{equation}{section}

\numberwithin{equation}{section}

\begin{document}

\title{Complete spectrum of the Robin eigenvalue problem on the ball
\thanks{Research supported by NNSF of China (No. 12371110).}}

\author{Guowei Dai\thanks{Corresponding author.
\newline
School of Mathematical Sciences, Dalian University of Technology, Dalian, 116024, P.R. China
\newline
\text{~~~~ E-mail}: daiguowei@dlut.edu.cn, sunyingxin2023@mail.dlut.edu.cn.}, Yingxin Sun}
\date{}
\maketitle

\renewcommand{\abstractname}{Abstract}

\begin{abstract}
We investigate the following Robin eigenvalue problem
\begin{equation*}
\left\{
\begin{array}{ll}
-\Delta u=\mu u\,\, &\text{in}\,\, B,\\
\partial_\texttt{n} u+\alpha u=0 &\text{on}\,\, \partial B
\end{array}
\right.
\end{equation*}
on the unit ball of $\mathbb{R}^N$. We obtain the complete spectral structure of this problem.
In particular, for $\alpha>0$, the first eigenvalue is $k_{\nu,1}^2$ and the second eigenvalue is $k_{\nu+1,1}^2$, where $k_{\nu+l,m}$ is the $m$th positive zero of $kJ_{\nu+l+1}(k)-(\alpha+l) J_{\nu+l}(k)$. Moreover, when $\alpha\in(-l,1-l)$ with any $l\in \mathbb{N}$, one has $l$ negative (strictly increasing) eigenvalues $-\widehat{k}_{\nu+i,1}^2$ with $i\in\{0,\ldots,l-1\}$ where $\widehat{k}_{\nu+l,1}$ denotes the unique zero of $\alpha I_{\nu+l}(k)+lI_{\nu+l}(k)+kI_{\nu+l+1}(k)$;
while, for $\alpha=-l$, besides $l$ negative (increasing) eigenvalues, $0$ is also an eigenvalue.
%In particular, when $\alpha\in(-1,0)$, the first eigenvalue is $-\widehat{k}_{\nu,1}^2$, and the second eigenvalue is exactly $k_{\nu+1,1}^2$.
%Furthermore, for $\alpha=-1$, the first eigenvalue is $-\widehat{k}_{\nu,1}^2$ and the second eigenvalue is exactly $0$.
\end{abstract}


\emph{AMS Subjection Classification(2020):} 34L20; 34L10; 35J05; 35P10; 35P20

%\tableofcontents
% ------------------------------------------------------------------------------------------------------------Introduction

\section{Introduction}
\bigskip
\quad\, Consider the following heat conduction equation
\begin{equation*}
w_t=\Delta w,\nonumber
\end{equation*}
where $w$ denotes the temperature as a function of position $x\in \Omega$ with $\Omega\subseteq \mathbb{R}^N$ and time $t$.
The radiation of heat from a homogeneous body $\Omega$ with boundary $\partial\Omega$ into an infinite medium maintained at zero temperature is described by the boundary condition
\begin{equation*}
\partial_\texttt{n} w+\sigma w=0,\nonumber
\end{equation*}
where $\sigma$ is a positive physical constant and $\texttt{n}$ is the unit outer normal to $\partial \Omega$. This condition indicates that the rate of change of temperature in the direction of the inner normal is proportional to the jump in temperature from the exterior to the interior of the body.
Using separation of variables for a given initial condition, one has the following Robin eigenvalue problem
\begin{equation}\label{Robin eigenvalue problem}
\left\{
\begin{array}{ll}
-\Delta u=\mu u\,\, &\text{in}\,\, \Omega,\\
(1-\beta)\partial_\texttt{n} u+\beta u=0 &\text{on}\,\, \partial \Omega,
\end{array}
\right.
\end{equation}
with $\beta\in \mathbb{R}$, which is also known as the elastically supported membrane problem.
This problem with positive boundary parameter is also related to the wave equation to elastically restoring boundaries.
Negative values of the Robin parameter appear in a model for surface superconductivity \cite{Giorgi}.
For $N=1$, it is well known that problem (\ref{Robin eigenvalue problem}) has and only has a sequence of simple eigenvalues,
and the eigenfunction corresponding to the $k$th eigenvalue possesses exactly $k-1$ simple zeros in $\Omega$ \cite{Coddington}.
For $N=2,3$ and when $\Omega$ is a ball, Courant and Hilbert \cite{Courant} outlined a method for calculating eigenvalues and eigenfunctions to this problem, though without providing specific calculation conclusion.
For general or special domains, this problem has been extensively studied by many mathematicians; we refer the reader to \cite{Henrot} and the references therein for a comprehensive overview.

When $\beta=1$, problem (\ref{Robin eigenvalue problem}) degenerates into the following fixed membrane problem
\begin{equation}\label{eigenvalue problem30}
\left\{
\begin{array}{ll}
-\Delta u=\lambda u\,\, &\text{in}\,\, \Omega,\\
u=0 &\text{on}\,\, \partial \Omega.
\end{array}
\right.
\end{equation}
Let $\lambda_n$ (with multiplicity) denote the $n$th eigenvalue of problem (\ref{eigenvalue problem30}) for $n \in \mathbb{N}$.
The spectral theory for the fixed membrane problem on balls or intervals is well-established; see, for instance, \cite{Courant} for dimensions $N \leq 3$, and \cite{Chavel} for $N \geq 4$.
For $N=2$, Payne, P\'{o}lya and Weinberger (PPW) \cite{PPW} established the inequality
\begin{equation*}
\lambda_{n+1}-\lambda_n\leq\frac{2}{n}\sum_{i=1}^n\lambda_i, \,\,\,n\in \mathbb{N}.
\end{equation*}
Hile and Protter \cite{HileProtter} later extended this result to arbitrary dimensions by proving that
\begin{equation*}
\sum_{i=1}^n\frac{\lambda_i}{\lambda_{n+1}-\lambda_i}\geq\frac{nN}{4}, \,\,\,n\in \mathbb{N}.
\end{equation*}
Yang \cite{Yang} obtained the following more sharp inequality
\begin{equation*}
\sum_{i=1}^n\left(\lambda_{n+1}-\lambda_i\right)\left(\lambda_{n+1}-\left(1+\frac{4}{N}\right)\lambda_i\right)\leq0, \,\,\,n\in \mathbb{N}.
\end{equation*}
In particular, for $N=2$ and $n=1$, all of the preceding inequalities imply that
\begin{equation*}
\frac{\lambda_2}{\lambda_1}\leq3.
\end{equation*}
PPW conjectured that the upper bound $3$ could be improved to the value attained when $\Omega$ is a disk, i.e.,
\begin{equation*}
\frac{\lambda_2}{\lambda_1}\leq \frac{\lambda_2}{\lambda_1}\Big|_{\text{disk}}=\frac{j_{1,1}^2}{j_{0,1}^2}\approx2.5387.
\end{equation*}
The PPW conjecture has been included in problem lists of Yau \cite[problem 77]{Yau}.
Yau commented that the truth of the conjecture means that the shape of a drum (whether it is circular or not) can be determined from its first two tones.
The upper bound of $3$ has been progressively improved by several mathematicians.
Specifically, Brands \cite{Brands} obtained $2.686$, de Vries \cite{Vries} improved it to $2.658$, and Chiti \cite{Chiti} further refined it to $2.586$.
For any dimension, the conjecture can be stated that
\begin{equation*}
\frac{\lambda_2}{\lambda_1}\leq \frac{\lambda_2}{\lambda_1}\Big|_{n-\text{ball}}=\frac{j_{n/2,1}^2}{j_{n/2-1,1}^2}.
\end{equation*}
This PPW conjecture was proved by Ashbaugh and Benguria \cite{AshbaughBenguria91, AshbaughBenguria92, AshbaughBenguria92CMP}.

For $\beta\neq1$, we write $\alpha=\beta/(1-\beta)\in \mathbb{R}$. Let $\mu_n$ (counting the multiplicity) with $n\in \mathbb{N}$ be the $n$th eigenvalue of the following elastically supported membrane problem
\begin{equation}\label{eigenvalue problem3}
\left\{
\begin{array}{ll}
-\Delta u=\mu u\,\, &\text{in}\,\, \Omega,\\
\partial_\texttt{n} u+\alpha u=0 &\text{on}\,\, \partial \Omega.
\end{array}
\right.
\end{equation}
In 2001, parallel to the fixed membrane problem, for $N=2$, Payne and Schaefer \cite{PayneSchaefer}
proposed the following conjecture.
\\ \\
\textbf{Payne-Schaefer conjecture.} \emph{The ratio
$\mu_2/\mu_1$ achieves its maximum for the disk for all values of $\alpha$ or for a range of
values of $\alpha$}.\\

\noindent They proved that $\mu_2/\mu_1\leq3$ for $\alpha\geq\alpha_*$ with some positive constant $\alpha_*$. In 2003, Henrot \cite[Open problem 15]{Henrot0} restated this conjecture by proposing the following open problem.
\\ \\
\textbf{Open problem.} \emph{For what values of $\alpha$ does the ratio
$\mu_2/\mu_1$ achieve its maximum for the disk?}\\

Although the original conjecture is formulated in two dimensions, it admits a natural extension to higher dimensions---namely, whether the ratio $\mu_2/\mu_1$ is maximized by the ball for all, or at least a range of, values of $\alpha$.

Recent years have seen active research on Robin eigenvalues for specific domains.
Freitas and Laugesen \cite{Freitas} investigated the first and second eigenvalues on balls (see also \cite{FreitasLaugesen} for the two-dimensional case).
Langford and Laugesen \cite{Langford} extended this analysis to geodesic disks in Euclidean, hyperbolic, and spherical spaces, using techniques that avoid special functions.
For positive boundary parameters $\alpha > 0$, Freitas \cite{Freitas2021} derived sharp estimates for the first eigenvalue on balls.
In the case of negative parameters, Antunes, Freitas, and Krej\v{c}i\v{r}\'{\i}k \cite{Antunes} and Freitas and Krej\v{c}i\v{r}\'{\i}k \cite{Freitas2015} established some bounds for $\mu_1$.
In \cite{Laugesen}, Laugesen studied various optimization problems for the first two eigenvalues on rectangular domains.

Unlike the fixed membrane problem, the ratio $\mu_2/\mu_1$ for a ball depends on the boundary parameter $\alpha$.
Building upon the aforementioned results, this work aims to systematically characterize the complete spectrum of the Robin eigenvalue problem on the ball.
A key objective is to derive the precise expression for $\mu_2/\mu_1$ on the ball, which is expected to provide a foundation for resolving the aforementioned conjecture and related open problems.

From now on, let $B$ denote the unit ball in $\mathbb{R}^N$ centered at the origin, and set $\Omega = B$.
Although our arguments extend to balls of arbitrary radius, we consider only the unit ball for simplicity.
To achieve our main objective, we present some more general results.
\\ \\
\textbf{Theorem 1.1.}
\emph{For problem (\ref{eigenvalue problem3}) with $\Omega=B$ and $N\geq2$, there exists a unique sequence of positive eigenvalues $\mu_{l,m}$ for $l\in\mathbb{N}\cup\{0\}$ and $m\in \mathbb{N}$ which has the following properties:}

(a) \emph{when $\alpha>-l$, the exact value of $\mu_{l,m}$ is $k_{\nu+l,m}^2$, where $k_{\nu+l,m}$ is the $m$th positive zero of $kJ_{\nu+l+1}(k)-(\alpha+l) J_{\nu+l}(k)$; when $\alpha=-l$, $\mu_{l,1}=0$ and $\mu_{l,m}$ is $k_{\nu+l,m}^2$ with $m\geq2$; when $\alpha<-l$, $\mu_{l,1}=-\widehat{k}_{\nu+l,1}^2$ where $\widehat{k}_{\nu+l,1}$ denotes the unique zero of $\alpha I_{\nu+l}(k)+lI_{\nu+l}(k)+kI_{\nu+l+1}(k)$, and $\mu_{l,m}$ is $k_{\nu+l,m}^2$ with $m\geq2$};

(b) \emph{when $\alpha>-l$ or $m\geq2$, the eigen-subspace corresponding to $\mu_{l,m}$ is generated by}
\begin{equation*}
u_{l,m}=r^{-\nu}J_{\nu+l}\left(k_{\nu+l,m}r\right)G_l(\xi);
\end{equation*}
\emph{for $\alpha=-l$, the eigen-subspace corresponding to $\mu_{l,1}=0$ is generated by $r^{-\alpha}G_{-\alpha}(\xi)$;
while, for $\alpha<-l$, the eigen-subspace corresponding to $\mu_{l,1}=-\widehat{k}_{\nu+l,1}^2<0$ is generated by $r^{-\nu}I_{\nu+l}\left(\widehat{k}_{\nu+l,1}r\right)G_{l}(\xi)$,
where $G_l(\xi)$ is any eigenfunction corresponding to $\kappa_l=l(l+N-2)$ on $N-1$ dimensional sphere $\mathbb{S}^{N-1}$;
in particular, $\mu_{0,1}$ is simple with positive radial eigenfunction, for any $m\geq2$,
$\mu_{0,m}$ is simple and the corresponding eigenfunction is radially symmetric and
has exactly $m-1$ simple zeros};

(c) \emph{for any fixed $l$, $\mu_{l,m}$ is strictly increasing with respect to $m$ and converges to infinity, meanwhile, for each fixed $l\in\mathbb{N}\cup\{0\}$, when $\alpha\geq-l $ and $m\geq2$, it is also strictly increasing with respect to $l$};

(d) \emph{when $\alpha\geq -l$ or $m\geq2$, one has that
\begin{equation*}
0\leq k_{\nu+l,1}<j_{\nu+l,1},\,\,j_{\nu+l+1,m-1}<k_{\nu+l,m}<j_{\nu+l,m};
\end{equation*}
for $\alpha>0$, the first eigenvalue is $k_{\nu,1}^2$ and the second eigenvalue is exactly $k_{\nu+1,1}^2$. Moreover, $k_{\nu,1}$ and $k_{\nu+1,1}$ are strictly increasing with respect to $\alpha$, and as $\alpha$ tends to $+\infty$, $k_{\nu,1}$ approaches $j_{\nu,1}$ and $k_{\nu+1,1}$ approaches $j_{\nu+1,1}$};

(e) \emph{when $\alpha\in(-l,1-l)$ with any $l\in \mathbb{N}$, one has $l$ negative (strictly increasing) eigenvalues $-\widehat{k}_{\nu+i,1}^2, i\in\{0,\ldots,l-1\}$;
while, for $\alpha=-l$ with any $l\in \mathbb{N}$, besides $l$ negative (increasing) eigenvalues, $0$ is also an eigenvalue.
In particular, when $\alpha\in(-1,0)$, the first eigenvalue is $-\widehat{k}_{\nu,1}^2$ and the second eigenvalue is exactly $k_{\nu+1,1}^2$; for $\alpha=-1$, the first eigenvalue is $-\widehat{k}_{\nu,1}^2$ and the second eigenvalue is exactly $0$}.
\\

Previous studies \cite{Antunes, Freitas, Freitas2015, Freitas2021, FreitasLaugesen, Langford, Laugesen} have established the existence and derived bounds for the first two eigenvalues. In this work, we go further by obtaining the full spectrum. Specifically, we provide explicit formulas (via Bessel functions) for the first two eigenvalues and determine the precise count of negative eigenvalues based on the boundary parameter. These results serve as a valuable supplement to the existing literature and provide a basis for future studies. It should be noted that some of our findings regarding the first two eigenvalues may overlap with the aforementioned references, yet we include them here for the sake of a self-contained presentation.

We employ the method of separation of variables (see, e.g., \cite{Courant, Freitas}) to establish the existence of eigenvalues and eigenfunctions. To further characterize their properties, we draw upon specific properties of Bessel functions. This approach enables us to derive explicit expressions for the first two eigenvalues. Specifically, applying Theorem 1.1 for $\alpha>0$ on the domain $B$ yields the ratio
\begin{equation*}
\frac{\mu_2}{\mu_1}\Big|_B=\frac{k_{\nu+1,1}^2}{k_{\nu,1}^2}.
\end{equation*}
For practical reference, we provide computed approximate values (computed using mathematical software) of $k_{\nu+l,1}$ and the resulting ratios $\mu_2/\mu_1$ for various $\alpha$ and $l$ in both $2$D and $3$D cases (see Table 1). Note that $\nu=N/2-1$.
\begin{table}[h]
\centering
\begin{tabular}{c|c|c|c|c|c|c|c|c}
               \hline
               $l$ & $\nu$ & $\alpha=1$ & $\alpha=2$ & $\alpha=3$ & $\alpha=4$ & $\alpha=5$ & $\alpha=100$ & $\alpha=1000$ \\ \hline
               $0$ & $0$          & 1.25578 & 1.59945 & 1.78866 & 1.90808 & 1.98981 & 2.38090 & 2.40242 \\
               $0$ & $1/2$ & 1.57080 & 2.02876 & 2.28893 & 2.45564 & 2.57043 & 3.11019 & 3.13845 \\
               $1$ & $0$          & 2.40483 & 2.73462 & 2.94960 & 3.09890 & 3.20752 & 3.79360 & 3.82788 \\
               $1$ & $$1/2$$ & 2.74371 & 3.14159 & 3.40561 & 3.59088 & 3.72638 & 4.44850 & 4.48892 \\
               \hline
               \multicolumn{2}{c|}{$\frac{\mu_2}{\mu_1}|_{2-\dim B}$} & 3.66726 & 2.92316 & 2.71938 & 2.63768 & 2.59846 & 2.53875 & 2.53874 \\
               \multicolumn{2}{c|}{$\frac{\mu_2}{\mu_1}|_{3-\dim B}$} & 3.05095 & 2.39794 & 2.21373 & 2.13832 & 2.10166 & 2.04575 & 2.04575 \\
               \hline
             \end{tabular}
             \caption{Approximate values of $k_{\nu+l,1}$ and $\mu_2/\mu_1$ for $l=0,1$ and $N=2,3$.}
\end{table}

As shown in Table 1, the values of $\mu_2/\mu_1$ become significantly smaller than the upper bound of $3$ established by Payne and Schaefer for general domains, particularly for large $\alpha$.
This observable dependence of the ratio on $\alpha$ validates the necessity of the constraints on $\alpha$ imposed in \cite{PayneSchaefer} to ensure $\mu_2/\mu_1 \leq 3$.
Consequently, our results provide a quantitative reference for refining such universal bounds.

For $\alpha=-1$, it follows directly that $\frac{\mu_2}{\mu_1}\Big|B=0$ when $N\geq2$.
For $\alpha\in (-1,0)$, the ratio is given by
\begin{equation*}
\frac{\mu_2}{\mu_1}\Big|_B=-\frac{k_{\nu+1,1}^2}{\widehat{k}_{\nu,1}^2}.
\end{equation*}
 Here we also provide approximate values (computed using mathematical software) of $k_{\nu+1,1}$, $\widehat{k}_{\nu,1}$ and the corresponding $\mu_2/\mu_1$ values for different values of $\alpha$ in the $2$D and $3$D cases (see Table 2).
\begin{table}[h]
\centering
\begin{tabular}{c|c|c|c|c|c|c}
               \hline
                & $\nu$ & $\alpha=-0.1$ & $\alpha=-0.3$ & $\alpha=-0.5$ & $\alpha=-0.7$ & $\alpha=-0.9$ \\
                \hline
               $k_{\nu+1,1}$ & $0$          & 1.76104 & 1.57883 & 1.35660 & 1.06842 & 0.62721  \\
               $k_{\nu+1,1}$ & $1/2$ & 1.98891 & 1.77934 & 1.52553 & 1.19873 & 0.70207 \\
               $\widehat{k}_{\nu,1}$ & $0$     & 0.45286 & 0.80454 & 1.06569 & 1.29403 & 1.50599 \\
               $\widehat{k}_{\nu,1}$ & $$1/2$$ & 0.55323 & 0.97767 & 1.28784 & 1.55477 & 1.79867 \\
               \hline
               \multicolumn{2}{c|}{$\frac{\mu_2}{\mu_1}|_{2-\dim B}$} & $-15.12204$ & $-3.85102$ & $-1.62047$ & $-0.68170$ & $-0.17345$\\
               \multicolumn{2}{c|}{$\frac{\mu_2}{\mu_1}|_{3-\dim B}$} & $-12.92465$ & $-3.31233$ & $-1.40319$ & $-0.59444$ & $-0.15236$ \\
               \hline
             \end{tabular}
             \caption{Approximate values of $k_{\nu+l,1}$, $\widehat{k}_{\nu,1}$ and $\mu_2/\mu_1$ for $N=2$ and $3$.}
\end{table}

In \cite[Conjecture B]{Laugesen}, Laugesen conjectured that $\mu_2/\mu_1$ is decreasing with respect to $\alpha>0$ on any bounded Lipschitz domain. Our numerical results, as presented in Tables 1 and 2, are consistent with this conjecture, at least within the ball domain.

\section{Spectrum of Robin problem in one dimension}

\bigskip
\quad\, To further study the spectrum of Robin problem on the ball, we first consider the one-dimensional case.
For $\Omega=(0,1)$, problem (\ref{eigenvalue problem3}) simplifies to
\begin{equation}\label{eigenvalue problem4}
\left\{
\begin{array}{ll}
-u''=\mu u,\,\, x\in(0,1),\\
-u'(0)+\alpha u(0)=0,\\
u'(1)+\alpha u(1)=0.
\end{array}
\right.
\end{equation}
It is well known \cite{Courant, Ince} that problem (\ref{eigenvalue problem4}) has and only has a sequence of simple eigenvalues $\mu_k$
with $\mu_1<\mu_2<\cdots\rightarrow+\infty$. For $\mu\geq0$, Bucur, Freitas and Kennedy \cite[Chapter 4]{Henrot} gave the following relation of eigenvalue $\mu$ and the parameter $\alpha$
\begin{equation*}
\alpha^2+2\alpha\sqrt{\mu}\cot\left(\sqrt{\mu}\right)-\mu=0
\end{equation*}
or the equivalent form
\begin{equation*}
\alpha_{\pm}=-\sqrt{\mu}\cot\left(\sqrt{\mu}\right)\pm\sqrt{\mu\csc^2\left(\sqrt{\mu}\right)}.
\end{equation*}
Meanwhile, for $\mu\leq0$, they also obtained that
\begin{equation*}
\alpha_+=-\sqrt{-\mu}\tanh\left(\frac{\sqrt{-\mu}}{2}\right)
\end{equation*}
and
\begin{equation*}
\alpha_-=-\sqrt{-\mu}\coth\left(\frac{\sqrt{-\mu}}{2}\right).
\end{equation*}
By analyzing the asymptotic behavior of $\alpha$ with respect to $\mu$, they obtained a rough image of the eigenvalue $\mu$ with respect to the parameter $\alpha$.
We pursue a more complete understanding by deriving their exact values, determining the number of negative eigenvalues, and obtaining the explicit expressions of the eigenfunctions.
\\ \\
\textbf{Proposition 2.1.}
\emph{For problem (\ref{eigenvalue problem3}) with $\Omega=(0,1)$, there exists a unique sequence of positive eigenvalues $\mu_m$ for $m\in \mathbb{N}$ which has the following properties:}

(a) \emph{$\mu_1>0$ when $\alpha>0$, $\mu_1=0$ for $\alpha=0$, $\mu_1<0$ and $\mu_2>0$ when $\alpha\in(-2,0)$, $\mu_1<0$ and $\mu_2=0$ for $\alpha=-2$, $\mu_2<0$ and $\mu_3>0$ when $\alpha<-2$. Moreover, one has that $(m-1)^2\pi^2<\mu_m<m^2\pi^2$ with $m\geq1$ for $\alpha>0$, $\mu_m=(m-1)^2\pi^2$ with $m\geq1$ for $\alpha=0$, $(m-2)^2\pi^2<\mu_m<(m-1)^2\pi^2$ with $m\geq2$ for $\alpha\in(-2,0)$ and $(m-2)^2\pi^2<\mu_m<(m-1)^2\pi^2$ with $m\geq3$ for $\alpha\leq-2$};

(b) \emph{when $\alpha\geq0$, the exact values of $\mu_m$ can be solved from}
\begin{equation*}
\alpha=-k\frac{\cos k}{\sin k}+ k\frac{1}{\vert \sin k\vert},
\end{equation*}
\emph{while, when $\alpha\leq0$, the exact values of $\mu_m$ can be solved from}
\begin{equation*}
\alpha=-k\frac{\cos k}{\sin k}- k\frac{1}{\vert \sin k\vert}
\end{equation*}
\emph{or}
\begin{equation*}
\alpha=-\sqrt{-\mu}\tanh\left(\frac{\sqrt{-\mu}}{2}\right)\,\,\text{or}\,\,\alpha=-\sqrt{-\mu}\coth\left(\frac{\sqrt{-\mu}}{2}\right);
\end{equation*}

(c) \emph{for $\alpha>0$ the eigen-subspace corresponding to $\mu_{m}$ is generated by}
 \begin{equation*}
u_m=\frac{\sqrt{\mu_{m}}}{\alpha}\cos\left(\sqrt{\mu_{m}} x\right)+\sin\left(\sqrt{\mu_{m}} x\right);
\end{equation*}
\emph{for $\alpha=0$ the eigen-subspace corresponding to $\mu_{m}$ is generated by $\cos((m-1)\pi x)$; for $\alpha\in(-2,0)$ the eigen-subspace corresponding to $\mu_{1}<0$ is generated by}
\begin{equation*}
\widetilde{u}_1=e^{-\sqrt{-\mu_1}x}-\frac{\alpha+\sqrt{-\mu_1}}{\alpha-\sqrt{-\mu_1}}e^{\sqrt{-\mu_1}x}
\end{equation*}
\emph{and the eigen-subspace corresponding to $\mu_{m}$ with $m\geq2$ is generated by $u_m$; for $\alpha=-2$ the eigen-subspace corresponding to $\mu_{1}<0$ is generated by the above $\widetilde{u}_1$,
corresponding to $\mu_{2}=0$ is generated by $1-2x$ and corresponding to $\mu_{m}$ with $m\geq3$ is generated by $u_m$; for $\alpha<-2$ the eigen-subspace corresponding to $\mu_{m}$ with $m=1,2$ is generated by}
\begin{equation*}
\widetilde{u}_m=e^{-\sqrt{-\mu_m}x}-\frac{\alpha+\sqrt{-\mu_m}}{\alpha-\sqrt{-\mu_m}}e^{\sqrt{-\mu_m}x}
\end{equation*}
\emph{and corresponding to $\mu_{m}$ with $m\geq3$ is generated by $u_m$}.
\\ \\
\textbf{Proof.} For $\mu \geq 0$, let $k = \sqrt{\mu}$. It is known that
\begin{equation*}
\alpha_{\pm}=-k\frac{\cos k}{\sin k}\pm k\frac{1}{\vert \sin k\vert}.
\end{equation*}
When $k\in(0,\pi)$, this simplifies to
\begin{equation*}
\alpha_{\pm}=-k\frac{\cos k}{\sin k}\pm k\frac{1}{\sin k}.
\end{equation*}
We now analyze the behavior of $\alpha_+$ on $(0, \pi)$. Noting that
\begin{equation*}
\lim_{k \to \pi^-} \frac{1 - \cos k}{\sin k} = +\infty \quad \text{and} \quad \lim_{k \to 0^+} \frac{1 - \cos k}{\sin k} = 0,
\end{equation*}
we obtain $\lim_{k \to \pi^-} \alpha_+ = +\infty$ and $\lim_{k \to 0^+} \alpha_+ = 0$. Furthermore, a direct computation gives
\begin{equation*}
\alpha_+'(k) = \frac{(1 - \cos k)(k + \sin k)}{\sin^2 k}, \quad k \in (0, \pi),
\end{equation*}
which implies that $\alpha_+$ is strictly increasing in $(0,\pi)$ which starts from $0$ and approaches to $+\infty$ as $k\rightarrow \pi^-$.

For $\alpha_{-}$ on $(0,\pi)$, we find
\begin{equation*}
\lim_{k\rightarrow \pi^-}\alpha_{-}(k)=-\pi\lim_{k\rightarrow \pi^-}\frac{1+\cos k}{\sin k}=0,
\end{equation*}
and
\begin{equation*}
\lim_{k\rightarrow 0^+}\alpha_{-}(k)=-\lim_{k\rightarrow 0^+}k\frac{1+\cos k}{\sin k}=-2.
\end{equation*}
Furthermore, the derivative
\begin{equation*}
\alpha_-'(k)=\frac{(1+\cos k)(k-\sin k)}{\sin^2 k}, \quad k\in(0,\pi),
\end{equation*}
is strictly positive, indicating that $\alpha_-$ is strictly increasing in $(0,\pi)$ which starts from $-2$ and joins to $0$.

By repeating the above analysis process, for each $m\in \mathbb{N}$, we find that $\alpha_{+}$ is strictly increasing in $(m\pi,(m+1)\pi)$, starts from $0$ and approaches $+\infty$ as $k\rightarrow ((m+1)\pi)^-$, and that $\alpha_{-}$ is also strictly increasing in $(m\pi,(m+1)\pi)$, starts from $-\infty$ and joins to $0$ as $k\rightarrow ((m+1)\pi)^-$.

For $\mu\leq0$, we have known that
\begin{equation*}
\alpha_+=-\sqrt{-\mu}\tanh\left(\frac{\sqrt{-\mu}}{2}\right)
\end{equation*}
and
\begin{equation*}
\alpha_-=-\sqrt{-\mu}\coth\left(\frac{\sqrt{-\mu}}{2}\right).
\end{equation*}
Differentiating $\alpha_+$ and $\alpha_-$ with respect to $\mu$ respectively, we have that
\begin{equation*}
  \alpha_+'=\frac{1}{2\sqrt{-\mu}}\left[\tanh\left(\frac{\sqrt{-\mu}}{2}\right)
  +\frac{\sqrt{-\mu}}{2}\text{sech}^2\left(\frac{\sqrt{-\mu}}{2}\right)\right]>0
\end{equation*}
and
\begin{equation*}
\begin{aligned}
  \alpha_-'&=\frac{1}{2\sqrt{-\mu}}\left[\coth\left(\frac{\sqrt{-\mu}}{2}\right)
  -\frac{\sqrt{-\mu}}{2}\text{csch}^2\left(\frac{\sqrt{-\mu}}{2}\right)\right]\\
  &=\frac{1}{4\sqrt{-\mu}}
  \left(\sinh\left(\sqrt{-\mu}\right)-\sqrt{-\mu}\right)\text{csch}^2\left(\frac{\sqrt{-\mu}}{2}\right)>0
  \end{aligned}
\end{equation*}
for $\mu<0$.
Thus we find that both $\alpha_+$ and $\alpha_-$ are strictly monotonically increasing, $\alpha_+>\alpha_-$ and $\alpha_+$ starts from $-\infty$ and joins to $0$ while $\alpha_-$ starts from $-\infty$ and joins to $-2$ in $(-\infty,0]$.

Based on the aforementioned properties of $\alpha_{\pm}$, we present a schematic diagram (see Figure 1; note that $k=\sqrt{\mu}\geq 0$ for $\mu\geq 0$ and $k=-\sqrt{-\mu}\leq0$ for $\mu\leq0$).
\begin{figure}[ht]
\centering
\includegraphics[width=0.6\textwidth]{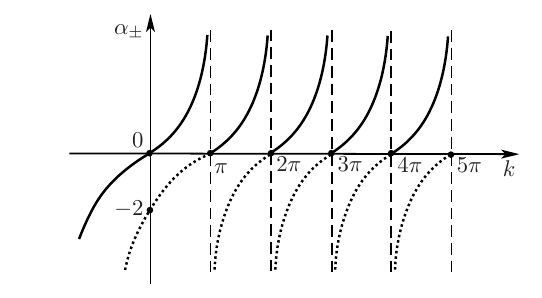}
\caption{The schematic diagram of $\alpha_+$ (solid lines) and $\alpha_-$ (dotted lines).}
\end{figure}
We see that
$\mu_1>0$ when $\alpha>0$, $\mu_1=0$ for $\alpha=0$, $\mu_1<0$ and $\mu_2>0$ when $\alpha\in(-2,0)$, $\mu_1<0$ and $\mu_2=0$ for $\alpha=-2$, $\mu_2<0$ and $\mu_3>0$ when $\alpha<-2$.
The range of positive eigenvalues $\mu_m$ is evident from the diagram, which establishes Property (a). Furthermore, Property (b) can be directly derived from the characteristics of $\alpha_{\pm}$, while Property (c) follows from the classical Euler undetermined function method.\qed
\\

As far as we know, the critical value $\alpha=-2$ represents a novel finding. This value marks the threshold at which the number of negative eigenvalues changes: when $\alpha<-2$, there exist exactly two negative eigenvalues, while for $\alpha\in[-2,0)$, there is exactly one.
From the monotonicity and asymptotic behavior of $\alpha_\pm$, we can deduce that both $\mu_1$ and $\mu_2$ tend to $-\infty$ as $\alpha\rightarrow-\infty$. In particular, for $\alpha=-3$, from Proposition 2.1 we obtain that $\mu_1\approx-10.52118$ and $\mu_2\approx-6.63412$.
For $\alpha>0$, we have that $\mu_1\in(0,\pi)$ and $\mu_2\in(\pi,2\pi)$. In this case, Freitas and J.B. Kennedy \cite{FreitasKennedy} obtained more refined upper and lower bound estimates for the first two eigenvalues.

\section{Spectrum of Robin problem on the ball}

\bigskip
\quad\, In this section, we establish the spectral structure of the Robin eigenvalue problem on the ball.
For $N=2$ and $\mu\geq0$, Bucur, Freitas and Kennedy \cite[Chapter 4]{Henrot} derived the following transcendental relation between $\alpha$ and $\mu$
\begin{equation*}
\alpha J_k\left(\sqrt{\mu}\right)+\frac{1}{2}\sqrt{\mu}\left[J_{k-1}\left(\sqrt{\mu}\right)-J_{k+1}\left(\sqrt{\mu}\right)\right]=0,
\end{equation*}
where $J_k$ denotes the Bessel function. For $\mu<0$, they also obtained that
\begin{equation*}
\alpha=-\frac{\sqrt{-\mu}\left[I_{k-1}\left(\sqrt{-\mu}\right)+I_{k+1}\left(\sqrt{-\mu}\right)\right]}{2I_k\left(\sqrt{-\mu}\right)},
\end{equation*}
where $I_k$ is the modified Bessel function.
Using the above relations, they further obtained the asymptotic expansion of $\mu$ or $\alpha$.

On the basis of the above classic conclusions, we further provide the exact values of eigenvalues on the ball with any $N\geq2$, as well as several properties of their corresponding eigenfunctions.
\\ \\
\textbf{Proof of Theorem 1.1.} Let $u$ be the eigenfunction corresponding to eigenvalue $\mu$.
If $u$ has the form $v(r)G(\xi)$ for $r\in[0,1]$ and $\xi\in \mathbb{S}^{N-1}$, from \cite{Chavel} we have that
\begin{equation*}
r^{1-N}\left(r^{N-1}v'\right)'G+r^{-2}v\Delta G+\mu vG=0
\end{equation*}
for $r\in(0,1)$.
It is well known that the distinct eigenvalues of the Laplace-Beltrami operator on $\mathbb{S}^{N-1}$ are given by $\kappa_l = l(l + N - 2)$ for $l \in \mathbb{N} \cup \{0\}$, and that the corresponding eigenfunctions are the $l$-th order spherical harmonics $G_l(\xi)$.
Thus, $G_l$ satisfies
\begin{equation*}
\Delta G_l + \kappa_l G_l = 0,
\end{equation*}
and the radial function $v$ obeys
\begin{equation*}
v'' + \frac{N-1}{r} v' + \left( \mu - \frac{\kappa_l}{r^2} \right) v = 0.
\end{equation*}

For $\mu>0$, define $k>0$ by $\mu=k^2$. Let $\tau=k r$ and $v(r)=z(\tau)$,
the above equation can be transformed into
\begin{equation*}
z''+\frac{N-1}{\tau}z'+\left(1-\frac{l(l+N-2)}{\tau^{2}}\right)z=0.
\end{equation*}
Let $J(\tau)=\tau^{\nu}z$. Then $J(\tau)$ satisfies the equation
\begin{equation*}
J''+\frac{1}{\tau}J'+\left(1-\frac{(\nu+l)^2}{\tau^{2}}\right)J=0,
\end{equation*}
the solution of which is $J(\tau)=J_{\nu+l}(\tau)$.
It follows that
\begin{equation*}
v(r)=r^{-\nu}J_{\nu+l}(kr),
\end{equation*}
up to a constant factor.
From the boundary condition we derive that $v'(1)+\alpha v(1)=0$, i.e.,
\begin{equation*}
\alpha J_{\nu+l}(k)+lJ_{\nu+l}(k)-kJ_{\nu+l+1}(k)=0
\end{equation*}
where we use the formula $J_\nu'(x)=-J_{\nu+1}(x)+{\nu}/{x}J_\nu(x)$ (see \cite{Olver}).

Since $J_{\nu+l}(k)$ and $J_{\nu+l+1}(k)$ cannot vanish simultaneously,
we have that
\begin{equation*}
\alpha +l-k\frac{J_{\nu+l+1}}{J_{\nu+l}}(k)=0.
\end{equation*}
This leads to
\begin{equation*}
\alpha=k\frac{J_{\nu+l+1}}{J_{\nu+l}}(k)-l:=\widetilde{h}_{\nu+l}(k).
\end{equation*}
The above equation also derived by Freitas \cite[Section 3]{Freitas2021} or \cite[Section 5]{Freitas} for the special case of $l=0$ and $\alpha>0$.
From the conclusions of \cite[Theorem B ]{Lorch} or \cite{Ismail}, we deduce that $\widetilde{h}_{\nu+l}(0)=-l$ and
$\widetilde{h}_{\nu+l}(k)$ is strictly increasing in $\left(0,j_{\nu+l,1}\right)$, it has poles $j_{\nu+l,m}$ and is also strictly increasing between any two adjacent poles.

We now consider the case of negative eigenvalues. For $\mu < 0$, let $t = -\mu > 0$. Applying separation of variables again yields, up to a constant factor, the radial solution
\begin{equation*}
v(r)=r^{-\nu}I_{\nu+l}\left(\sqrt{t}r\right).
\end{equation*}
Substituting this into the boundary condition leads to
\begin{equation*}
\alpha I_{\nu+l}(k)+kI_{\nu+l}'(k)-\nu I_{\nu+l}(k)=0,
\end{equation*}
where $k=\sqrt{t}$.
Using the identity for the derivative of the modified Bessel function \cite{Olver}
\begin{equation*}
I_\nu'(x)=I_{\nu+1}(x)+\frac{\nu}{x}I_\nu(x),
\end{equation*}
we obtain that
\begin{equation*}
\alpha I_{\nu+l}(k)+lI_{\nu+l}(k)+kI_{\nu+l+1}(k)=0.
\end{equation*}
It follows that
\begin{equation*}
\alpha=-k\frac{I_{\nu+l+1}}{I_{\nu+l}}(k)-l:=\widehat{h}_{\nu+l}(k).
\end{equation*}
From the results in \cite[Theorem C]{Lorch}, we deduce that $\widehat{h}_{\nu+l}(0)=-l$ and that $\widehat{h}_{\nu+l}(k)$ is strictly decreasing in $\left(0,+\infty\right)$.

Finally, for $\mu = 0$, a similar analysis shows that the radial equation becomes
\begin{equation*}
v''+\frac{N-1}{r}v'-\frac{\kappa_l}{r^{2}}v=0.
\end{equation*}
Using the Euler undetermined function method, we find that
\begin{equation*}
\begin{aligned}
v(r)=r^{l},
\end{aligned}
\end{equation*}
up to a constant factor. In view of the boundary condition, one has that $l=-\alpha$, which also have been obtained in \cite[Section 5]{Freitas}.

In conclusion, we define
\begin{equation*}
{h}_{\nu+l}(k)=\left\{
\begin{array}{ll}
\widetilde{h}_{\nu+l}(k)\,\, &\text{for}\,\, \mu>0,\\
-l  &\text{for}\,\, \mu=0,\\
\widehat{h}_{\nu+l}(k)\,\, &\text{for}\,\, \mu<0.
\end{array}
\right.
\end{equation*}
Hence, from the above properties of ${h}_{\nu+l}(k)$, we derive that, for each fixed $l\in \mathbb{N}\cup\{0\}$, ${h}_{\nu+l}(k)=\alpha$ has and only has a sequence of roots $k_{\nu+l,m}$ for any $m\in \mathbb{N}$ (which depends on the position of $\alpha$,
see following Figure 2).
\begin{figure}[htbp]
\centering
\subfloat[$l=0$ and $\mu>0$]{
\includegraphics[scale=0.85]{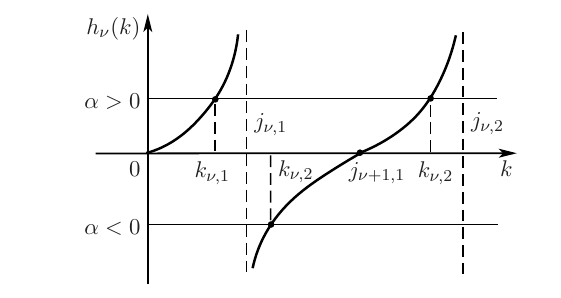}
}
\hfill
\subfloat[$l=0$ and $\mu<0$]{
\includegraphics[scale=0.85]{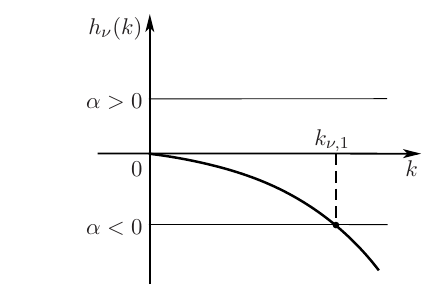}
}
\quad
\subfloat[$l\neq0$ and $\mu>0$]{
\includegraphics[scale=0.85]{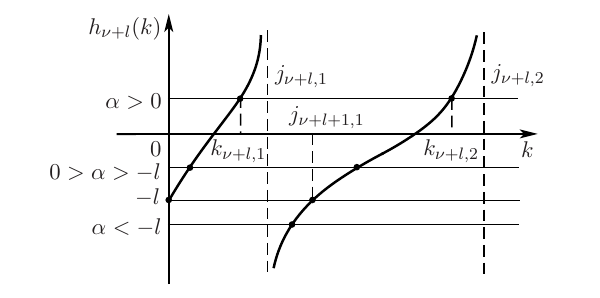}
}
\hfill
\subfloat[$l\neq0$ and $\mu<0$]{
\includegraphics[scale=0.85]{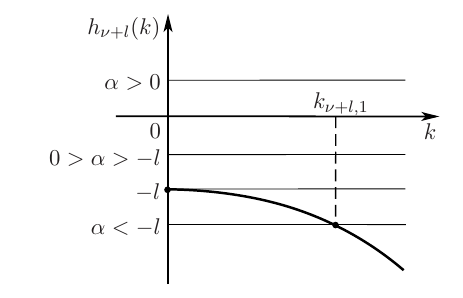}
}
\caption{The schematic diagram of $h_{\nu+l}(k)$.}
\end{figure}
For $\alpha>-l$, $\mu_{l,m}=k_{\nu+l,m}^2$ are eigenvalues of problem (\ref{eigenvalue problem3}) where $k_{\nu+l,m}=\widetilde{k}_{\nu+l,m}$ denotes the $m$th zero of $\alpha=\widetilde{h}_{\nu+l}(k)$.
For $\alpha=-l$, $\mu_{l,1}=0$ and $\mu_{l,m}=k_{\nu+l,m}^2>0$ with $m\geq2$.
For $\alpha<-l$, $\mu_{l,1}=-\widehat{k}_{\nu+l,1}^2<0$ and $\mu_{l,m}=k_{\nu+l,m}^2>0$ with $m\geq2$ where $k_{\nu+l,1}=\widehat{k}_{\nu+l,1}$ denotes the unique zero of $\alpha=\widehat{h}_{\nu+l}(k)$.
This establishes Property (a).

When $\alpha>-l$ or $m\geq2$, using  $\mu_{l,m}=k_{\nu+l,m}^2$, we obtain that
\begin{equation*}
v(r)=r^{-\nu}J_{\nu+l}\left(k_{\nu+l,m}r\right):=v_{l,m}.
\end{equation*}
Then
\begin{equation*}
u_{l,m}:=v_{l,m}(r)G_l(\xi)
\end{equation*}
is the eigenfunction corresponding to $\mu_{l,m}$.
For $\alpha=-l$, the eigen-subspace corresponding to $\mu_{l,1}=0$ is generated by $r^{-\alpha}G_{-\alpha}(\xi)$.
While, for $\alpha<-l$, the eigen-subspace corresponding to $\mu_{l,1}=-\widehat{k}_{\nu+l,1}^2<0$ is generated by $r^{-\nu}I_{\nu+l}\left(\widehat{k}_{\nu+l,1}r\right)G_{l}(\xi)$,
%We have shown that the eigen-subspace corresponding to $\mu_{l,m}$ is generated by
% \begin{equation*}
%u_{l,m}=r^{-\nu}J_{\nu+l}\left(k_{\nu+l,m}r\right)G_l(\xi),
%\end{equation*}
where $G_l(\xi)$ is any eigenfunction corresponding to $\kappa_l=l(l+N-2)$ on $N-1$ dimensional sphere $\mathbb{S}^{N-1}$.

Similar to that of \cite{Chavel} with obvious changes we know that the function-space, $L$, consisting of the span, in $L^2(B)$, of all eigenfunctions of $B$ obtained by the above procedure, is dense in $L^2(B)$.
Consequently, the eigenvalue sequence ${\mu_{l,m}}$ constitutes the full spectrum of problem (\ref{eigenvalue problem3}).

In particular, $\mu_{0,1}$ is simple with positive radial eigenfunction.
For any $m\geq2$, the eigen-subspace corresponding to $\mu_{0,m}$ is one-dimensional and generated by $v_{0,m}(r)$.
Hence, $\mu_{0,m}$ is simple and the corresponding eigenfunction $v_{0,m}(r)$ is radially symmetric.
It has been known (see, for example, \cite{Chavel} or \cite[Theorem 2.1 of Chap. 8]{Coddington}) that $v_{0,m}(r)$ has exactly $m-1$ simple zeros.
This completes the proof of property (b).

For any fixed $l\in \mathbb{N}\cup\{0\}$ and $\alpha\geq-l$, from the properties of $h_{\nu+l}(k)$, we derive that
the sequence $k_{\nu+l,m}$ is strictly increasing with respect to $m$, and hence so is $\mu_{l,m} = k_{\nu+l,m}^2$. Furthermore, the interlacing inequalities
\begin{equation}\label{interlacing inequalities}
0 \leq k_{\nu+l,1} < j_{\nu+l,1}, \quad j_{\nu+l,m-1}<k_{\nu+l,m}<j_{\nu+l,m} \,\,\text{for}\,\,m\geq2.
\end{equation}
together with $j_{\nu+l,m} \to +\infty$ as $m \to \infty$, imply that $k_{\nu+l,m} \to +\infty$ as $m \to \infty$. Consequently, $\mu_{l,m} \to +\infty$ as $m \to \infty$.

Moreover, from \cite[Theorem A]{Lorch} we know that $\widetilde{h}_{\nu+l}(k)$ is strictly decreasing with respect to $l$. Combining this with the monotonicity of $h_{\nu+l}(k)$ with respect to $k$ implies that $k_{\nu+l,m}$ with $\alpha\geq-l$ or $m\geq2$ is strictly increasing with respect to $l$. Thus, when $\alpha\geq-l$ or $m\geq2$, $\mu_{l,m}$ is strictly increasing with respect to $l$, which verifies property (c).

We finally prove property (d). When $\alpha\geq-l$ or $m\geq2$, inequalities (\ref{interlacing inequalities}) hold.
Furthermore, when $\alpha>0$, the identity $h_{\nu+l}(j_{\nu+l+1,m-1}) = -l$ leads to
\begin{equation*}
j_{\nu+l+1,m-1}<k_{\nu+l,m}<j_{\nu+l,m}
\end{equation*}
for $m\geq2$.
When $\alpha>0$, from the above properties one sees that
\begin{equation*}
k_{\nu+1,1}<j_{\nu+1,1}<k_{\nu,2},
\end{equation*}
and hence we obtain that
\begin{equation*}
\mu_1=k_{\nu,1}^2>0
\end{equation*}
and
\begin{equation*}
\mu_2=k_{\nu+1,1}^2.
\end{equation*}
Moreover, it is clear from the monotonicity of $h_{\nu+l}$ that $k_{\nu,1}$ and $k_{\nu+1,1}$ are strictly increasing with respect to $\alpha$. As $\alpha$ tends to $+\infty$, $k_{\nu,1}$ approaches $j_{\nu,1}$ and $k_{\nu+1,1}$ approaches $j_{\nu+1,1}$. Thus we obtain the desired conclusions of (d).

Finally we show property (e). From the properties of ${h}_{\nu+l}(k)$ we can derive the following conclusions for $\alpha<0$.
When $\alpha\in(-1,0)$, the first eigenvalue is $-\widehat{k}_{\nu,1}^2$ where $\widehat{k}_{\nu,1}$ denotes the unique zero of $\alpha I_{\nu}(k)+kI_{\nu+1}(k)$, and the second eigenvalue is exactly $k_{\nu+1,1}^2$ where $k_{\nu+1,1}$ is the first positive zero of $kJ_{\nu+2}(k)-(\alpha+1) J_{\nu+1}(k)$.
For $\alpha=-1$, the first eigenvalue is $-\widehat{k}_{\nu,1}^2$ where $\widehat{k}_{\nu,1}$ denotes the unique zero of $-I_{\nu}(k)+kI_{\nu+1}(k)$, and the second eigenvalue is exactly $0$.
In general, for any given $l\in \mathbb{N}$ with $l\geq 2$, when $\alpha\in(-l,1-l)$, the first eigenvalue is $\min_{i\in\{0,\ldots,l-1\}}\left\{-\widehat{k}_{\nu+i,1}^2\right\}$
where $\widehat{k}_{\nu+i,1}$ denotes the unique zero of $\alpha I_{\nu+i}(k)+iI_{\nu+i}(k)+kI_{\nu+i+1}(k)$.
While, for $l\in \mathbb{N}$ with $l\geq 2$ and $\alpha=-l$, the first eigenvalue is $\min_{i\in\{0,\ldots,l-1\}}\left\{-\widehat{k}_{\nu+i,1}^2\right\}$
where $\widehat{k}_{\nu+i,1}$ denotes the unique zero of $-lI_{\nu+i}(k)+iI_{\nu+i}(k)+kI_{\nu+i+1}(k)$.
From \cite[Lemma 10]{Freitas} we derive that $\widehat{h}_{\nu+l}(k)$ is strictly decreasing with respect to $l$.
This implies that $\widehat{k}_{\nu+i,1}^2$ is strictly decreasing with respect to $i$.
Hence, $-\widehat{k}_{\nu+i,1}^2$ is strictly increasing with respect to $i$.
In particular, the first eigenvalue is just $-\widehat{k}_{\nu,1}^2$ and the second eigenvalue is $-\widehat{k}_{\nu+1,1}^2$.
%In particular, when $\alpha\in(-1,0)$, the first eigenvalue is $-\widetilde{k}_{\nu,1}^2$ where $\widetilde{k}_{\nu,1}$ denotes the unique zero of $\alpha I_{\nu}(k)+kI_{\nu+1}(k)$, and the second eigenvalue is exactly $k_{\nu+1,1}^2$.
%Furthermore, for $\alpha=-1$, the first eigenvalue is $-\widetilde{k}_{\nu,1}^2$ where $\widetilde{k}_{\nu,1}$ denotes the unique positive zero of $-I_{\nu}(k)+kI_{\nu+1}(k)$, and the second eigenvalue is exactly $0$.
\qed\\

Theorem 1.1 reveals a fundamental difference from the one-dimensional case (see Proposition 2.1): as $\alpha\rightarrow-\infty$, the problem admits infinitely many negative eigenvalues. Since $I_{\nu+1}(k)/I_\nu(k)$ increases monotonically to $1$ as $k \rightarrow +\infty$, the function $\widehat{h}_{\nu+l}(k)$ decreases monotonically to $-\infty$ in the same limit. This leads to $\mu_{l,1}\rightarrow-\infty$ as $\alpha \rightarrow -\infty$. In particular, the first eigenvalue decreases monotonically to $-\infty$.

Anyway, we have obtained the complete spectral structure of problem (\ref{eigenvalue problem3}) by utilizing the properties of Bessel functions. In particular, we obtain the exact values of eigenvalues and the exact expressions for the basis of the eigen-subspace. These conclusions themselves are also interesting.
\\ \\
\textbf{The conflicts of interest statement and Data Availability statement.}
\bigskip\\
\indent There is not any conflict of interest.
Data sharing not applicable to this article as no datasets were generated or analysed during the current study.
%\\ \\
%\textbf{Acknowledgment}
%\bigskip\\
%\indent The author expresses deep gratitude to Professor Mark S. Ashbaugh for raising the monotonicity problem of $k_\nu(a)$ and for his guidance and discussion on this question.

\bibliographystyle{amsplain}
\makeatletter
\def\@biblabel#1{#1.~}
\makeatother

%\bibliography{mybib2}

\providecommand{\bysame}{\leavevmode\hbox to3em{\hrulefill}\thinspace}
\providecommand{\MR}{\relax\ifhmode\unskip\space\fi MR }
\providecommand{\MRhref}[2]{%
  \href{http://www.ams.org/mathscinet-getitem?mr=#1}{#2}
}
\providecommand{\href}[2]{#2}

\end{document}